\documentclass[11pt]{article}
\newcommand{\bbR}{{\rm I\hspace{-0.7mm}R}}

\usepackage{latexsym}

\begin{document}

\centerline{\bf \large Hausdorff Dimension For Level Sets And $k$-Multiple Times }

\vskip .1in

\centerline{\tt By MING YANG}

\vskip .1in

\indent{\footnotesize
We compute the Hausdorff dimension of the zero-set of an additive L\'evy process.
As an application, we obtain
the Hausdorff dimension formula for the set of $k$-multiple times of a L\'evy process.}

\vskip .1in
\noindent
{\small \noindent 2000 {\it Mathematics Subject Classification}. Primary
60G60, 60G51; secondary 60G17.
  
\noindent
{\it Key words and phrases}. Additive L\'evy processes, Hausdorff dimension, level sets, intersection times,
$k$-multiple times.}

\vskip .15in
\centerline{\textsc{1. Introduction}}
\vskip .15in

Let
$X^1_{t_1},~X^2_{t_2},\cdots,X^N_{t_N}$ be $N$ independent L\'evy
processes in $\bbR^d$ with their respective L\'evy exponents 
$\Psi_j,~j=1,2,\cdots,N$. The random field
$$X_t=X^1_{t_1}+X^2_{t_2}+\cdots+X^N_{t_N},~~~~~t=(t_1,t_2,\cdots,t_N)\in\bbR^N_+$$
is called the additive
L\'evy process. Let $\lambda_d$ denote Lebesgue measure in $\bbR^d.$
Define
$E_1+E_2=\{x+y:x\in E_1,~y\in E_2\}$ 
for any two sets $E_1,~E_2$ of $\bbR^d.$ 
The following theorem
has recently been proved.

\vskip 0.15in
\noindent
{\bf Theorem 1.1}~~
{\it Let $X$ be any additive L\'evy process in $\bbR^d$ with
L\'evy exponent $(\Psi_1,\cdots,\Psi_N).$ Then
for any $F\in\mathcal{B}(\bbR^d)\backslash\{\emptyset\},$}
$$E\{\lambda_d(X(\bbR^N_+)+F)\}>0\Longleftrightarrow
\int_{\bbR^d}|\hat{\mu}(\xi)|^2
\prod_{j=1}^N\mbox{Re}\left(\frac{1}{1+\Psi_j(\xi)}\right)
d\xi<\infty\eqno(1.1)$$
{\it for some probability measure $\mu$ on $F$, where
$\hat{\mu}(\xi)=\int_{\bbR^d}e^{i\xi\cdot x}\mu(dx),~\xi\in\bbR^d.$} 

\vskip .15in
 
\noindent
Khoshnevisan, Xiao and Zhong [1] proved (1.1) under a sector condition. 
Yang [3] recently removed the sector condition. A special case of Theorem 1.1 is the following
theorem.

\vskip 0.15in
\noindent
{\bf Theorem 1.2}~~
{\it
Let $X$ be any additive L\'evy process in $\bbR^d$ with
L\'evy exponent $(\Psi_1,\cdots,\Psi_N).$
Then}
$$E\{\lambda_d(X(\bbR^N_+))\}>0\Longleftrightarrow
\int_{\bbR^d}
\prod_{j=1}^N\mbox{Re}\left(\frac{1}{1+\Psi_j(\xi)}\right)
d\xi<\infty.\eqno(1.2)$$
[Note that (1.1) still holds when $\bbR_+$ is replaced by 
$(0,\infty)$ for some or all
components of $\bbR^N_+.$]

\vskip .15in
 
\noindent
Theorem 1.2 is the key for the results of this paper. Meanwhile Theorem 1.1 is also needed in order to
relate symmetric stable L\'evy processes to stable subordinators.
We use
a subordination technique 
to solve the level set problem. To be more
precise, we give the formula for Hausdorff dimension of the level set.
Subsequently, the intersection time problem and in particular
the $k$-multiple time problem are solved. 

Finally, we mention a $q$-potential density criterion:
Let $X$ be an additive L\'evy process and assume that
$X$ has an a.e. positive $q$-potential density on $\bbR^d.$ Then
$$P\left\{
F\bigcap X((0,\infty)^N)\neq\emptyset\right\}>0\Longleftrightarrow
E\left\{\lambda_d(F-X((0,\infty)^N))\right\}>0.\eqno(1.3)$$
The argument is elementary but crucially hinges on the
property: 
$X_{b+t}-X_b,~t\in\bbR^N_+$ (independent of $X_b$) can be replaced by $X$ 
for all $b\in\bbR^N_+$; moreover, the second condition ``a.e. positive on $\bbR^d$"
is absolutely necessary for the direction $\Longleftarrow$ in (1.3); 
see for example
Proposition 6.2 of [1].

\vskip .15in
\centerline{\textsc{2. Level Sets and Intersection Times}}

\vskip .15in

Let $Z$ be an $N$-parameter additive L\'evy process in $\bbR^d.$ The
set $Z^{-1}(0)=\{t\in(0,\infty)^N:Z_t=0\}$ is called the level set 
(at $0\in\bbR^d$; otherwise better known as the zero-set).
The level set problem, in short, is to compute $\dim_HZ^{-1}(0).$

First,  
by (1.3) 
and Theorem 1.2 we obtain immediately that

\vskip 0.15in
\noindent
{\bf Theorem 2.1}~~
{\it Let $\{Z;~\Psi_1,\cdots,\Psi_N\}$ be an $N$-parameter additive L\'evy process in $\bbR^d.$ 
Assume that $Z$ has an a.e. positive $q$-potential density on $\bbR^d$ for some $q\ge0.$
Then}
$$P(Z^{-1}(0)\neq\emptyset)>0
\Longleftrightarrow
\int_{\bbR^d}
\prod_{j=1}^N\mbox{Re}\left(\frac{1}{1+\Psi_j(\xi)}\right)
d\xi<\infty.\eqno(2.1)$$

\vskip .15in
Define for $\alpha\in(0,1)$
the following $(N,N)$ additive L\'evy process: 
$$\sigma^{\alpha}_t=\sigma^1_{t_1}+\sigma^2_{t_2}+\cdots+\sigma^N_{t_N},$$  
where
$$\sigma^j_{t_j}=(\sigma^{j,1}_{t_j},\sigma^{j,2}_{t_j},
\cdots,\sigma^{j,N}_{t_j}),~~~j=1,2,\cdots,N$$
and the $\sigma^{j,l}_{t_j}, ~1\le j\le N,~1\le l\le N,$ are $N^2$
i.i.d. standard $\alpha$-stable subordinators with the Laplace exponent
$\lambda^{\alpha}.$ We call $\sigma^{\alpha}_t$ the saturated additive 
$\alpha$-stable subordinator. We note that for each $t\in(0,\infty)^N$,
$\sigma^{\alpha}_t$ has a density positive everywhere on $(0,\infty)^N$.   
Let $\Psi^{\alpha}_j$ be the L\'evy exponent
of $\sigma^j_{t_j}.$ 
A standard $\alpha$-stable subordinator with the Laplace exponent
$\lambda^{\alpha}$ has L\'evy exponent 
$$|\lambda|^{\alpha}-i|\lambda|^{\alpha}\mbox{sgn}(\lambda)\theta,~~~\lambda\in\bbR,$$
where $\theta=\tan(\alpha\frac{\pi}{2})\in(0,\infty).$
Thus,
$$\Psi^{\alpha}_j(\xi)=
\sum^N_{l=1}|\xi_l|^{\alpha}-i\left(\sum^N_{l=1}
|\xi_l|^{\alpha}\mbox{sgn}(\xi_l)\right)\theta,~~~\xi=(\xi_1,\cdots,\xi_N)\in\bbR^N.$$

\vskip 0.15in
\noindent
{\bf Lemma 2.2}~~
{\it Let $X_t$ be a L\'evy process in $\bbR^d$ with L\'evy exponent
$\Psi$ and let $\sigma^1,~\sigma^2,\cdots,\sigma^n$ be $n$ independent subordinators also
independent of $X$. Then $X_{\sigma^1_{t_1}+\sigma^2_{t_2}+\cdots+\sigma^n_{t_n}}$
is an $n$-parameter additive L\'evy process with L\'evy exponent
$(\Psi^{\sigma^1},\Psi^{\sigma^2},\cdots,\Psi^{\sigma^n}),$ where
$\Psi^{\sigma^j}$ denotes the 
$\sigma^j$-subordination of $\Psi.$}

\vskip 0.15in
\noindent
{\bf Proof}~~By independence,
$$Ee^{i\xi\cdot X_{\sigma^1_{t_1}+\sigma^2_{t_2}+\cdots+\sigma^n_{t_n}}}
=
Ee^{-(\sigma^1_{t_1}+\sigma^2_{t_2}+\cdots+\sigma^n_{t_n})\Psi(\xi)}
=
\prod^n_{j=1}Ee^{-\sigma^j_{t_j}\Psi(\xi)}
=\prod^n_{j=1}e^{-t_j\Psi^{\sigma^j}(\xi)}.\eqno\Box$$

\vskip .15in
The next lemma follows readily from Lemma 2.2 
and the subordination formula.

\vskip 0.15in
\noindent
{\bf Lemma 2.3}~~
{\it $Z\circ{\sigma}^{\alpha}$ is an $N$-parameter additive L\'evy process in $\bbR^d$
(a sum of $N$ i.i.d. L\'evy processes)
with L\'evy exponent} 
$$\left\{\sum^N_{j=1}[\Psi_j(\xi)]^{\alpha},~\sum^N_{j=1}[\Psi_j(\xi)]^{\alpha},\cdots,
~\sum^N_{j=1}[\Psi_j(\xi)]^{\alpha}\right\}.\eqno(2.2)$$

\vskip 0.15in
\noindent
{\bf Lemma 2.4}~~
{\it
For any Borel set $F\subset(0,\infty)^N,$}
$$P\left\{
F\bigcap \sigma^{\alpha}((0,\infty)^N)\neq\emptyset\right\}>0\Longleftrightarrow
E\left\{\lambda_N(F-\sigma^{\alpha}((0,\infty)^N))\right\}>0.\eqno(2.3)$$

\vskip 0.15in
\noindent
{\bf Proof}~~The direction 
$\Longrightarrow$
is easy since $\sigma^{\alpha}_t$ has density. We establish the
direction $\Longleftarrow,$ which is of some interest.
Assume that
$$E\left\{\lambda_N(F-\sigma^{\alpha}((0,\infty)^N))\right\}>0.$$
We extend $\sigma^{\alpha}$ to a truly saturated subordinator over
the entire time parameter space $\bbR^N.$ Due to the special property
of the subordinator, which maps $\bbR_+$ into $\bbR_+$ and
thanks to the way we chose to define $\sigma^{\alpha}$, this extension
is exactly the one it ought to be.
For example, in the second quadrant,
$$\bar{\sigma}^j=(-\bar{\sigma}^{j,1}_{-t_j},\bar{\sigma}^{j,2}_{t_j},
\cdots,\bar{\sigma}^{j,N}_{t_j}),~~~j=1,2,\cdots,N$$
where $t_j\in\bbR_+$, $\bar{\sigma}^{j,1}_{-t_j}\stackrel{d}{=}\sigma^{j,1}_{t_j},$
and the $\bar{\sigma}$ are i.i.d. copies of the $N^2$ $\sigma$'s in the first
quadrant. Then we have a saturated subordinator 
$\bar{\sigma}^{\alpha}=\bar{\sigma}^1+\cdots+\bar{\sigma}^N$ in the second quadrant.
After all quadrants have a saturated subordinator,
we see that the $2^N$ saturated subordinators are independent, disjoint in the interior
and each is a duplicate of $\sigma^{\alpha}$. Thus, we have defined a
process $\Sigma$ on $\bbR^N.$ Since on each open quadrant, the corresponing saturated subordinator
has a positive density everywhere for every $t$ in the same open quadrant, we see that for any $q>0,$
$\Sigma$ has 
an a.e. positive $q$-potential density on $\bbR^N.$
It is also clear that $\Sigma$ 
is a process of the property: 
$\Sigma_{b+t}-\Sigma_b,~t\in\bbR^N_+$ (independent of $\Sigma_b$)
can be replaced by
$\sigma^{\alpha}$
for all $b\in\bbR^N$.
Since $F\subset(0,\infty)^N,$ we can also make each quadrant have
a copy of $F$ through reflection. Let $F^*$ be the union of
the $2^N$ disjoint copies of $F$. Then clearly,
$$P\left\{
F\bigcap \sigma^{\alpha}((0,\infty)^N)\neq\emptyset\right\}>0\Longleftrightarrow
P\left\{
F^*\bigcap \Sigma(\bbR^N_I)\neq\emptyset\right\}>0,\eqno(2.4)$$
where $\bbR^N_I$ is the union of the interiors of the quadrants of $\bbR^N.$
What we have done is the so-called ``coloring the coordinate axes of the
first quadrant".
Since $F\subset F^*$,
$$E\left\{\lambda_N(F^*-\sigma^{\alpha}((0,\infty)^N))\right\}>0.$$
Since $\Sigma$ has positive $q$-potential density a.e. on $\bbR^N,$
and more importantly since $\Sigma$
has the property: $\Sigma_{b+t}-\Sigma_b,~t\in\bbR^N_+$ (independent of $\Sigma_b$)
can be replaced by
$\sigma^{\alpha}$
for all $b\in\bbR^N$,
applying the standard $q$-potential density argument shows that
$$P\left\{
F^*\bigcap \Sigma(\bbR^N_I)\neq\emptyset\right\}>0.\eqno\Box$$

\vskip .15in

Let 
$$S^{\alpha}_t=S^1_{t_1}+S^2_{t_2}+\cdots+S^N_{t_N}$$  
be the standard 
$N$-parameter additive $\alpha$-stable
L\'evy process in $\bbR^N$ for $\alpha\in(0,1);$
that is,
the $S^j$ are independent standard $\alpha$-stable
L\'evy processes in $\bbR^N$
with the common L\'evy exponent
$|\xi|^{\alpha}.$
A non-trivial result in multi-parameter potential theory is given by the last lemma.

\vskip 0.15in
\noindent
{\bf Lemma 2.5}~~
{\it 
For any Borel set $F\subset(0,\infty)^N,$ }  
$$P\left\{
F\bigcap S^{\alpha}((0,\infty)^N)\neq\emptyset\right\}>0\Longleftrightarrow
P\left\{
F\bigcap \sigma^{\alpha}((0,\infty)^N)\neq\emptyset\right\}>0.$$

\vskip 0.15in
\noindent
{\bf Proof}~~ By Lemma 2.4,
$$P\left\{
F\bigcap \sigma^{\alpha}((0,\infty)^N)\neq\emptyset\right\}>0\Longleftrightarrow
E\left\{\lambda_N(F-\sigma^{\alpha}((0,\infty)^N))\right\}>0.\eqno(2.5)$$
Since $S^{\alpha}$ has an a.e. positive $1$-potential density on $\bbR^d,$
by (1.3) 
$$P\left\{
F\bigcap S^{\alpha}((0,\infty)^N)\neq\emptyset\right\}>0\Longleftrightarrow
E\left\{\lambda_N(F-S^{\alpha}((0,\infty)^N))\right\}>0.\eqno(2.6)$$

Let $\Psi^{\alpha}_j$ be the L\'evy exponent
of $\sigma^j_{t_j}.$ 
We show that
there are two constants $A_1,~A_2\in(0,\infty)$ such that
$$A_1\left(\frac{1}{1+|\xi|^{\alpha}}\right)^N
\le
\prod_{j=1}^N\mbox{Re}\left(\frac{1}{1+\Psi_j^{\alpha}(\xi)}\right)
\le A_2
\left(\frac{1}{1+|\xi|^{\alpha}}\right)^N,~~\forall~\xi\in\bbR^N.\eqno(2.7)$$
We have seen that
$$\Psi^{\alpha}_j(\xi)=
\sum^N_{l=1}|\xi_l|^{\alpha}-i\left(\sum^N_{l=1}
|\xi_l|^{\alpha}\mbox{sgn}(\xi_l)\right)\theta,~~~\xi=(\xi_1,\cdots,\xi_N)\in\bbR^N,$$
where $\theta=\tan(\alpha\frac{\pi}{2})\in(0,\infty).$
It follows from simple calculations that
$$\left(\frac{1}{1+\theta^2}\right)\cdot
\frac{1}{1+\sum^N_{l=1}|\xi_l|^{\alpha}}
<\mbox{Re}\left(\frac{1}{1+\Psi_j^{\alpha}(\xi)}\right)
<\frac{1}{1+\sum^N_{l=1}|\xi_l|^{\alpha}}.$$
Now, (2.7) follows from the well-known two-sided inequality
$$1+|\xi|^{\alpha}\le
1+\sum^N_{l=1}|\xi_l|^{\alpha}
\le
C_{\alpha,N}
(1+|\xi|^{\alpha})$$
where 
$C_{\alpha,N}\in(1,\infty)$ is a constant depending only on $\alpha$ and $N.$
Thanks to (2.7), Theorem 1.1 implies that (considering the probability measures on $-F$),
$$E\left\{\lambda_N(F-S^{\alpha}((0,\infty)^N))\right\}>0\Longleftrightarrow
E\left\{\lambda_N(F-\sigma^{\alpha}((0,\infty)^N))\right\}>0.\eqno\Box$$

\vskip .15in 

 For each $\beta\in(0,N)$, choose a  
$\sigma^{1-\beta/N}$ independent of $Z$ and define
$Z\circ{\sigma}^{1-\beta/N}.$
The following theorem is the main result of this paper.

\vskip 0.15in
\noindent
{\bf Theorem 2.6}~~
{\it
Let $\{Z;~\Psi_1,\cdots,\Psi_N\}$ be an $N$-parameter 
additive L\'evy process in $\bbR^d.$
Assume that for each $\beta\in(0,N)$, 
$Z\circ{\sigma}^{1-\beta/N}$ 
has
an a.e. positive $q$-potential density on $\bbR^d$ for some $q\ge0.$ ($q$ might depend on $\beta$.)
[A special case
is that if
for each $t\in(0,\infty)^N,$
$Z_t$ has an a.e. positive density on $\bbR^d,$ then 
$Z\circ{\sigma}^{1-\beta/N}$ has an a.e. positive $1$-potential 
density on $\bbR^d$ for all $\beta\in(0,N).$]  
If
$P(Z^{-1}(0)\neq\emptyset)>0$, then almost surely 
$\dim_HZ^{-1}(0)$ is a constant  
on $\{Z^{-1}(0)\neq\emptyset\}$ and} 
$$\dim_HZ^{-1}(0)=
\sup\left\{\beta\in(0,N):
\int_{\bbR^d}\left[
\mbox{Re}\left(\frac{1}{1+\sum^N_{j=1}[\Psi_j(\xi)]^{1-\beta/N}}\right)\right]^N
d\xi<\infty\right\}.\eqno(2.8)$$

\vskip 0.15in
\noindent
{\bf Proof}~~
According to the argument, Eq. (4.96)-(4.102), in {\it Proof of Theorem 3.2.} of
Khoshnevisan, Shieh, and Xiao [2], it suffices to show that
for all $\beta\in(0,N)$ and $S^{1-\beta/N}$ independent of $Z$,
$$P\left\{
Z^{-1}(0)\bigcap S^{1-\beta/N}((0,\infty)^N)\neq\emptyset\right\}>0\Longleftrightarrow
\int_{\bbR^d}\left[
\mbox{Re}\left(\frac{1}{1+\sum^N_{j=1}[\Psi_j(\xi)]^{1-\beta/N}}\right)\right]^N
d\xi<\infty.\eqno(2.9)$$
By Lemma 2.5, 
$$P\left\{
Z^{-1}(0)\bigcap S^{1-\beta/N}((0,\infty)^N)\neq\emptyset\right\}>0\Longleftrightarrow
P\left\{
Z^{-1}(0)\bigcap \sigma^{1-\beta/N}((0,\infty)^N)\neq\emptyset\right\}>0$$
$$\Longleftrightarrow
P\left\{0\in Z\circ{\sigma}^{1-\beta/N}((0,\infty)^N)\right\}>0.$$
Since $Z\circ{\sigma}^{1-\beta/N}$ 
has
an a.e. positive $q$-potential density,
by (1.3), Lemma 2.3, and Theorem 1.2
$$P\left\{0\in Z\circ{\sigma}^{1-\beta/N}((0,\infty)^N)\right\}>0
\Longleftrightarrow
\int_{\bbR^d}\left[
\mbox{Re}\left(\frac{1}{1+\sum^N_{j=1}[\Psi_j(\xi)]^{1-\beta/N}}\right)\right]^N
d\xi<\infty.\eqno\Box$$

\vskip .15in

Let 
$(X^1;~\Psi_1),$
$\cdots,~(X^k;~\Psi_k)$ 
be $k$ independent L\'evy processes in $\bbR^d$ for $k\ge 2$.
The set
$$T_k=\{(t_1,\cdots,t_k)\in(0,\infty)^k:X^1_{t_1}=\cdots=X_{t_k}^k\}$$
is called the intersection time set.
Define
$$Z_t=(X^2_{t_2}-X^1_{t_1},\cdots,X^k_{t_k}-X^{k-1}_{t_{k-1}}).$$
$Z$ is a $k$-parameter additive L\'evy process taking values
in $\bbR^{d(k-1)}.$
Clearly,
$$T_k=Z^{-1}(0).$$
Note 
that for any $\bbR^d$-valued random variable $X$ and $\xi_1,~\xi_2\in\bbR^d,$
$e^{i[(\xi_1,\xi_2)\cdot(X,-X)]}
=e^{i(\xi_1-\xi_2)\cdot X}.$ In particular,
the L\'evy process $(X^j,-X^j)$ has L\'evy exponent
$\Psi_j(\xi_1-\xi_2).$ It follows that
$Z^j$,
$1\le j\le k,$ has L\'evy exponent $\Psi_j(\xi_j-\xi_{j-1})$
with $\xi_0=\xi_k=0.$ 
By specializing Theorems 2.1 and 2.6, we have

\vskip 0.15in
\noindent
{\bf Theorem 2.7}~~
{\it
Let $(X^1;~\Psi_1),~\cdots,~(X^k;~\Psi_k)$ be $k$ independent L\'evy processes in $\bbR^d$ for
$k\ge 2$. 
Assume that $Z$ has an a.e. positive $q$-potential density for some $q\ge0.$ [A special case
is that if 
for each $j=1,\cdots,k,$
$X^j$ has a one-potential density $u_j^1>0,$ $\lambda_d$-a.e., then
$Z$ has an a.e. positive $1$-potential density on $\bbR^{d(k-1)}.$]  
Then}
$$P(T_k\neq\emptyset)>0
\Longleftrightarrow
\int_{\bbR^{d(k-1)}}
\prod_{j=1}^k\mbox{Re}\left(\frac{1}{1+\Psi_j(\xi_j-\xi_{j-1})}\right)
d\xi_1\cdots d\xi_{k-1}<\infty\eqno(2.10)$$
{\it with $\xi_0=\xi_k=0.$

Assume that for each $\beta\in(0,k)$, 
$Z\circ{\sigma}^{1-\beta/k}$ 
has
an a.e. positive $q$-potential density on $\bbR^{d(k-1)}$ for some $q\ge0,$
where $\sigma^{1-\beta/k}$ has $k$ parameters. ($q$ might depend on $\beta$.)
[A special case
is that if
for each $t\in(0,\infty)$ and each $j=1,\cdots,k,$
$X^j_t$ has 
an a.e. positive density on $\bbR^d,$ then 
$Z\circ{\sigma}^{1-\beta/k}$ has an a.e. positive $1$-potential 
density on $\bbR^{d(k-1)}$ for all $\beta\in(0,k).$]  
If
$P(T_k\neq\emptyset)>0$, then almost surely 
$\dim_HT_k$ is a constant  
on $\{T_k\neq\emptyset\}$ and} 
$$\dim_HT_k=
\sup\left\{\beta\in(0,k):
\int_{\bbR^{d(k-1)}}
\left[\mbox{Re}\left(\frac{1}{1+\sum^k_{j=1}
[\Psi_j(\xi_j-\xi_{j-1})]^{1-\beta/k}}\right)\right]^k
d\xi_1\cdots d\xi_{k-1}<\infty\right\}\eqno(2.11)$$
{\it with $\xi_0=\xi_k=0.$}

\vskip .15in

Let $X$ be a L\'evy process in $\bbR^d$ with L\'evy exponent $\Psi.$
The set for $k\ge2$ 
$$L_k=\{(t_1,\cdots,t_k)\in\bbR^k_+:~t_1,\cdots,t_k~~\mbox{are distinct,}~~
X_{t_1}=\cdots=X_{t_k}\}$$
is called the $k$-multiple time set of $X$.
It is a standard fact that
$L_k$ can be identified with
$T_k$  
as long as we replace the 
$X^j$
by the i.i.d. copies of $X.$ Thus, we have 
found the solution to the multiple-time problem:

\vskip 0.15in
\noindent
{\bf Theorem 2.8}~~
{\it Let $(X,~\Psi)$ be any L\'evy process in $\bbR^d$.
Assume that for each $\beta\in(0,k)$, 
$Z\circ{\sigma}^{1-\beta/k}$ 
has
an a.e. positive $q$-potential density on $\bbR^{d(k-1)}$ for some $q\ge0.$ ($q$ might depend on $\beta$.)
[A special case
is that if
for each $t\in(0,\infty),$
$X_t$ has 
an a.e. positive density on $\bbR^d,$ then 
$Z\circ{\sigma}^{1-\beta/k}$ has an a.e. positive $1$-potential 
density on $\bbR^{d(k-1)}$ for all $\beta\in(0,k).$]
Let
$L_k$ be the $k$-multiple-time set of $X$ for $k\ge 2$.  
If
$P(L_k\neq\emptyset)>0$, then almost surely
$\dim_HL_k$ is a constant  
on $\{L_k\neq\emptyset\}$ and}  
$$\dim_HL_k=
\sup\left\{\beta\in(0,k):
\int_{\bbR^{d(k-1)}}\left[\mbox{Re}\left(\frac{1}{1+\sum^k_{j=1}
[\Psi(\xi_j-\xi_{j-1})]^{1-\beta/k}}\right)\right]^k
d\xi_1\cdots d\xi_{k-1}<\infty\right\}\eqno(2.12)$$
{\it with $\xi_0=\xi_k=0.$}

\vskip .3in

\centerline{REFERENCES}
\vskip .2in

\noindent {[1]} Khoshnevisan, D., Xiao, Y. and Zhong, Y. (2003). 
Measuring the range of an additive 

\indent \indent
L\'evy process. {\it Ann. Probab.} {\bf 31}, 1097-1141.

\noindent {[2]} Khoshnevisan, D., Sheih, N.-R. and Xiao, Y. (2006).
Hausdorff dimension of the contours 

\indent \indent
of symmetric additive proceeses. {\it Probab. Th. Rel. Fields}. (To appear.)

\noindent {[3]} Yang, M. (2007). On a general theorem for additive L\'evy processes.
{\it Submitted.}

\end{document}